\documentclass[11pt]{article}
\usepackage{amsmath}
\usepackage{mathrsfs}
\usepackage{amsfonts}

\usepackage{amsfonts, amsmath, amssymb}
\usepackage{amssymb,amsfonts,amsmath,
latexsym, epsfig,cite, psfrag,eepic,colordvi}
\usepackage{amscd,graphics}
\usepackage{lineno}

\newcommand{\qed}{\hfill $\Box $}
\newcommand{\pf}{\noindent {\bf Proof.} }
\newcommand{\rem}{\noindent {\bf Remark} }

\newtheorem{theorem}{Theorem}[section]

\newtheorem{coro}[theorem]{Corollary}

\begin{document}

\title{A Note On Edge Connectivity and Parity Factor}

\author{Hongliang Lu\thanks{Corresponding email: luhongliang215@sina.com},
\\ {\small  School of Mathematics and Statistic}
\\ {\small Xi'an Jiaotong University, Xi'an 710049, PR China}
}

\date{}

\maketitle

\date{}

\maketitle

\begin{abstract}
In this paper, we investigate some parity factors by using
Lov\'asz's $(g,f)$-parity theorem. Let $m>0$ be an integer. Firstly,
we obtain a sufficient and necessary condition for some graphs to
have a parity factor with restricted minimum degree. Secondly, we
obtain some sufficient conditions for a graph to have a parity
factor with  minimum degree $m$ in term of edge connectivity.

\end{abstract}

\section{Introduction}

Let $G = (V,E)$ be a graph with vertex set $V (G)$ and edge set
$E(G)$. The number of vertices of a graph $G$ is called the
\emph{order} of $G$.  For a vertex $v$ of  graph $G$, the number of
edges of $G$ incident with $v$ is called the \emph{degree} of $v$ in
$G$ and is denoted by $d_{G}(v)$. For two subsets $S,T\subseteq
V(G)$, let $e_{G}(S,T)$ denote the number of edges of $G$ joining
$S$ to $T$. The \emph{edge cut} of $G$ is a subsets of edges whose
total removal renders the graph disconnected and the
\emph{edge-connectivity} is the size of a smallest edge cut.

An \emph{even (odd) factor} of $G$ is a spanning subgraph of $G$ in
which every vertex has even (odd, resp.) positive degree. Let $g, f
: V\rightarrow Z^+$ be two integer-valued function such that
$g(v)\leq f (v)$ and $g(v)\equiv f (v) \pmod 2$  for every $v\in V$.
Then a spanning subgraph $F$ of $G$ is called a $(g, f
)$-parity-factor, if $g(v)\leq d_F(v)\leq f (v)$ and $d_F(v)\equiv f
(v) \pmod 2$ for all $v\in V$. By the definitions, it is easy to see
that an  even (or odd) factor is also a special parity factor.


Petersen obtained a well-known result on the existence of a
2-factor.
\begin{theorem}[Petersen,\cite{Peter1}]\label{Peter}
Every bridgeless cubic graph has a 2-factor.
\end{theorem}

Fleischner extended Theorem \ref{Peter} to the existence of an even
factor.
\begin{theorem}[Fleischner,\cite{Fleischner}; Lov\'asz, \cite{Lovasz-Com}]\label{Fleischner}
If $G$ is a bridgeless graph with $\delta(G)\geq 3$, then $G$ has an
even factor.
\end{theorem}

For $(g,f)$-parity factor, Lov\'asz obtained a sufficient and
necessary condition.
\begin{theorem}[Lov\'asz \cite{Lovasz}]\label{Lovasz}
Let $G$ be a graph and let $g,f:V(G)\rightarrow N$ such that $g\leq
f$ and $g(v)\equiv f(v)\pmod 2$ for all $v\in V(G)$. Then $G$
contains a $(g,f)$-parity factor if and only if for any two disjoint
subsets $S$ and $T$,
\begin{align*}
 \eta(S,T)=g(T)-\sum_{x\in T}d_{G-S}(x)-f(S) +\tau(S,T)>0,
\end{align*}
where $\tau(S,T)$ denotes the number of components $C$, called
$g$-odd components of $G-(S\cup T)$ such that
$e_{G}(V(C),T)+f(V(C))\equiv 1 \pmod 2$.
\end{theorem}

In this paper, we shall study some parity factors by using
Lov\'asz's $(g,f)$-parity theorem.  
Firstly, for parity factor with minimum degree restricted, we obtain
a simple sufficient and necessary condition. Secondly, we generalize
Theorem \ref{Fleischner} and obtain some sufficient conditions for a
graph to have an even (or odd) factor in term of edge connectivity.

\section{Main Results}

\begin{theorem}\label{Main}
Let $G$ be a graph and $g:V(G)\rightarrow N$ be an integer-valued
function. Then $G$ contains a factor $F$ such that $d_F(v)\geq g(v)$
and $d_F(v)\equiv g(v)$ for all $v\in V(G)$ if and only if for any
subset $T\subseteq V(G)$,
\begin{align*}
g(T)-\sum_{x\in T}d_{G}(x) +\tau(T)\leq 0,
\end{align*}
where $\tau(T)$ denotes the number of components $C$, called $g$-odd
components of $G-T$ such that $e_{G}(V(C),T)+g(V(C))\equiv 1 \pmod
2$.
\end{theorem}

\pf Firstly, we prove necessity. Suppose that necessity does not
hold. Then there exists $T\subseteq V(G)$ such that
\begin{align*}
g(T)-\sum_{x\in T}d_{G}(x) +\tau(T)>0,
\end{align*}
where $\tau(T)$ denotes the number of components $C$, called $g$-odd
components of $G-T$ such that $e_{G}(V(C),T)+g(V(C))\equiv 1 \pmod
2$.  Let $F$ be a  parity factor of $G$ such that $d_F(v)\equiv
g(v)\pmod 2$ and $d_F(v)\geq g(v)$. By parity, $F$ misses at least
an edge from every $g$-odd component to $T$. Then we have
\begin{align*}
 \sum_{x\in T}d_{G}(x)-\tau(T)\geq \sum_{x\in T}d_{F}(x)\geq g(T),
\end{align*}
a contradiction.

Now we prove the sufficiency.  Let $f:V(G)\rightarrow N$ is  an
integer-valued function such that $f(v)\geq \Delta(G)+1$ and
$f(v)\equiv g(v)\pmod 2$ for all $v\in V(G)$. Suppose that the
sufficiency does not hold. Then $G$ contains no $(g,f)$-parity
factor. By Theorem \ref{Lovasz}, there exists two  disjoint subsets
$S$ and $T$ such that
\begin{align}\label{Parity_eq1}
  g(T)-\sum_{x\in T}d_{G-S}(x)-f(S) +\tau(S,T)>0,
\end{align}
where $\tau(S,T)$ denotes the number of components $C$, called
$g$-odd components of $G-(S\cup T)$ such that
$e_{G}(V(C),T)+g(V(C))\equiv 1 \pmod 2$. We choose $S, T$ such that
$S$ is minimal. Now we claim that  $S=\emptyset$. Otherwise, suppose
that $S\neq \emptyset$. Let $v\in S$ and $S'=S-v$. Then we have
\begin{align*}
\eta(S',T)&=g(T)-\sum_{x\in T}d_{G-S'}(x)-f(S') +\tau(S',T)\\
&\geq g(T)-\sum_{x\in
T}d_{G-S}(x)+e_G(v,T)-f(S)+f(v)+\tau(S,T)-(d_G(v)-e_G(v,T)+1)\\
&=g(T)-\sum_{x\in
T}d_{G-S}(x)-f(S)+\tau(S,T)+f(v)+2e_G(v,T)-d_G(v)-1,\\
&\geq g(T)-\sum_{x\in T}d_{G-S}(x)-f(S)+\tau(S,T)>0,
\end{align*}
contradicting to the minimality of $S$.

By (\ref{Parity_eq1}), then the result is followed. This completes
the proof. \qed

\begin{coro}\label{Main_even}
Let $G$ be a graph and $m>0$ be even. Then $G$ contains an even
factor $F$ such that $d_F(v)\geq m$ if and only if for any  subset
$T$,
\begin{align*}
m|T|-\sum_{x\in T}d_{G}(x) +\tau(T)\leq 0,
\end{align*}
where $\tau(T)$ denotes the number of components $C$, called $m$-odd
components of $G-T$ such that $e_{G}(V(C),T)\equiv 1 \pmod 2$.
\end{coro}

\begin{coro}\label{Main_odd}
Let $G$ be a graph and $m>0$ be odd. Then $G$ contains an odd factor
$F$ such that $d_F(v)\geq m$ if and only if for any  subset $T$,
\begin{align*}
m|T|-\sum_{x\in T}d_{G}(x) +\tau(T)\leq 0,
\end{align*}
where $\tau(T)$ denotes the number of components $C$, called $m$-odd
components of $G-T$ such that $e_{G}(V(C),T)+|V(C)|\equiv 1 \pmod
2$.
\end{coro}

\begin{theorem}\label{Theorem_even}
Let $m$ be even.  Let $G$ be an $m$-edge-connected graph with
minimum degree $m+1$. Then $G$ contains an even factor $F$ such that
$d_F(v)\geq m$.
\end{theorem}

\pf Suppose that the result does not hold. By Corollary
\ref{Main_even}, there exists $T\subseteq V(G)$ such that
\begin{align}\label{Even_eq1}
m|T|-\sum_{x\in T}d_{G}(x) +\tau(T)>0,
\end{align}
where $\tau(T)$ denotes the number of components $C$, called $m$-odd
components of $G-T$ such that $e_{G}(V(C),T)\equiv 1 \pmod 2$. Let
$C_1,\ldots,C_{\tau}$ denote these $m$-odd components of $G-T$. We
write $W=C_1\cup \cdots \cup C_{\tau}$.
 Now we discuss two cases.

\medskip
 {\it Case 1.} $\tau(T) \leq |T|$.
\medskip

Since $\delta(G)\geq m+1$, then we have
\begin{align*}
\sum_{x\in T}d_{G}(x)\geq (m+1)|T|.
\end{align*}
Hence
\begin{align}
m|T|-\sum_{x\in T}d_{G}(x) +\tau(T)\leq \tau(T)-|T|\leq 0,
\end{align}
contradicting to inequality (\ref{Even_eq1}).

\medskip
 {\it Case 2.} $\tau(T) > |T|$.
\medskip

Since $G$ is $m$-edge connected, then we have $e_G(V(C),T) \geq m$.
Note that $e_G(V(C),T)\equiv 1\pmod 2$ and $m$ is even. Hence
$e_G(V(C),T) \geq m+1$. So we have
\begin{align}
\sum_{x\in T}d_{G}(x)\geq e_G(T,V(W))\geq (m+1)\tau(T),
\end{align}
which implies
\begin{align}
m|T|-\sum_{x\in T}d_{G}(x) +\tau(T)\leq m|T|- (m+1)\tau(T)
+\tau(T)\leq 0,
\end{align}
a contradiction again. This completes this proof. \qed

\rem~\textbf{1:} The bound of Theorem \ref{Theorem_even} is tight.
Let $K_{2m}$ denote the complete graph with order $2m$. Take $m+1$
copies of $K_{2m}$ and add $m+1$ new vertices $v_1,\ldots,v_{m+1}$.
We write $T=\{v_1,\ldots,v_{m+1}\}$. Choose one vertex from each
copy of $K_{2m}$ and connect the vertex to every vertex of set $T$.
This results an $(m+1)$-edge-connected graph. Since
\begin{align}
(m+2)|T|-\sum_{x\in T}d_{G}(x) +\tau(T)= (m+2)|T|- (m+1)\tau(T)
+\tau(T)>0,
\end{align}
by Corollary \ref{Main_even}, then $G$ contains no an even factor
$F$ such that $d_F(v)\geq m+2$.

\begin{theorem}\label{Theorem_odd}
Let $m$ be odd.  Let $G$ be an $(m+1)$-edge-connected graph. Then
$G$ contains an odd factor $F$ such that $d_F(v)\geq m$.
\end{theorem}

\pf Suppose that the result does not hold. By Corollary
\ref{Main_odd}, there exists a subset $T$ such that
\begin{align}\label{Odd_eq1}
m|T|-\sum_{x\in T}d_{G}(x) +\tau(T)>0,
\end{align}
where $\tau(T)$ denotes the number of components $C$, called $m$-odd
components of $G-T$ such that $e_{G}(V(C),T)+|V(C)|\equiv 1 \pmod
2$.  We write $\tau=\tau(T)$. Let $C_1,\ldots,C_{\tau}$ denote these
$m$-odd components of $G-T$ and let $W=C_1\cup \cdots \cup
C_{\tau}$.

Since $G$ is $m+1$-edge-connected, then we have
\begin{align*}
\sum_{x\in T}d_{G}(x)\geq \max\{(m+1)|T|,(m+1)\tau\},
\end{align*}
which implies
\begin{align}
m|T|-\sum_{x\in T}d_{G}(x) +\tau\leq 0
\end{align}
a contradiction. This completes the proof. \qed

By the discussion in Remark 1, it is easy to show that the bound of
Theorem \ref{Theorem_odd} is also tight.

\end{document}